\documentclass[11pt]{article}
\usepackage{graphicx}
\usepackage{amsfonts}
\usepackage{bbm}
\usepackage{mathrsfs}
 \usepackage{mathrsfs,amsmath,amssymb,amsthm}
\usepackage[dvips]{color}
 \usepackage{amssymb}
 \usepackage{amsbsy}
 \usepackage[dvips]{color}
\allowdisplaybreaks
 \theoremstyle{plain}
\theoremstyle{remark}  \newtheorem{remark}{\noindent\mbox{Remark}}
 \theoremstyle{remark}\newtheorem{example}{\noindent\mbox{Example}}
 \theoremstyle{plain}\newtheorem{lemma}{\noindent\mbox{Lemma}}
\theoremstyle{plain} \newtheorem{theorem}{\noindent\mbox{Theorem}}
 \theoremstyle{plain}
 \theoremstyle{plain}
\theoremstyle{definition}

 \def\bq{\begin{equation}}
 \def\eq{\end{equation}}
 \def\eqn{\end{eqnarray}}
 \def\bqn{\begin{eqnarray}}
 \def\proof{\noindent{\it Proof.~~}}
 \def\qed{\hfill$\Box$\medskip}
 \def\rto{\rightarrow\infty}
 \def\z{\left}
 \def\y{\right}
 \def\no{\nonumber}

\begin{document}
 \title{\textbf{Coalescence for supercritical Galton-Watson processes with immigration}}                  
\author{  Hua-Ming \uppercase{Wang}$^{a},$ Lulu \uppercase{Li}$^b$ \& Huizi \uppercase{Yao}$^c$ \\\\
     \mbox{}$^{a,b,c}$School of Mathematics and Statistics, Anhui Normal University,\\
      Wuhu 241003, People's Republic of China\\
    E-mail: $^a$hmking@ahnu.edu.cn; $^b$lilulu91@163.com; $^c$yaohuiziyao@163.com }
\date{}
\maketitle%

\vspace{-.8cm}

\begin{center}
\begin{minipage}[c]{12cm}
\begin{center}\textbf{Abstract}\quad \end{center}

In this paper, we consider  Galton-Watson processes with immigration. Pick $i(\ge2)$ individuals randomly without replacement from the $n$-th generation and trace their lines of descent back in time till they coalesce into $1$ individual in a certain generation, which we denote by $X_{i,1}^n$ and is called the coalescence time. Firstly, we give the probability distribution of $X_{i,1}^n$ in terms of the probability generating functions of both the offspring distribution and the immigration law. Then by studying the limit behaviors of various functionals of the Galton-Watson process with immigration, we find the limit distribution of $X_{2,1}^n$ as $n\rto.$
\vspace{0.2cm}

\textbf{Keywords:}\  Galton-Watson processes; immigration; coalescence time.

\textbf{MSC 2010}: 60J80; 62E15
\vspace{0.2cm}

\end{minipage}
\end{center}

\section{Introduction}

Recent years, coalescence becomes an interesting research object in the community of branching processes.
 Athreya \cite{a12} studied the distribution of the coalescence time and its limit behavior of critical and subcritical Galton-Watson  processes whereas the supercritical case was considered by the same author in \cite{a12b}.
 Lambert \cite{lama} investigated the limit distribution of the coalescence time in the subcritical case for  some more general settings, including both discrete and continuous time and state space. Furthermore, Lambert and Popovic \cite{lp} defined a coalescent point
process, for which the coalescence time of two successive individuals alive at the same time is defined as
the first point mass in it, and study its limit.
For coalescence in Bellman-Harris and multitype branching processes, we refer the reader to   Athreya and Hong \cite{ah} and  Hong \cite{ha}. We also note that Grosjean and Huillet \cite{gh}  studied
a more general coalescence for Galton-Watson processes, that is,  coalescence of $i(\ge2)$ individuals  into $j(<i)$ individuals. For more results and details of  coalescence in branching processes, we refer the reader to Athreya \cite{a15} and references therein for a survey.

In this paper, we consider Galton-Watson branching processes with {\it immigration}.
To give the precise model, let $M$ be a $\mathbb Z^+:=\{0,1,2,...\}$-valued random variable with probability generating function (p.g.f.)  $f(z)=E(z^M),$ $z\in[0,1],$ which  serves as the offspring distribution and let $I$ be also a $\mathbb Z^+$-valued random variable with p.g.f. $g(z)=E(z^I),$ $z\in[0,1],$  which plays the role as the immigration law.
Set $N_0=0$ and for  $n \ge 0,$ define
\begin{equation}\label{nr}
  N_{n+1}=\sum_{i=1}^{N_n+I_n}\xi_{n,i},
\end{equation}
where $I_n,$ $n\geq 0$ are mutually independent and all distributed as $I,$ for $n\ge0,$ $I_n$ are independent of $\{N_n, N_{n-1},...,N_0\},$ and $\xi_{n,i},$ $n\geq 0,$ $i\geq 1$ are mutually independent and have the same distribution as $M.$ We mention that in (\ref{nr}) and the remainder of the paper, empty sum equals $0.$
We call the process $\{N_n\}_{n\ge0}$ a Galton-Watson  process with immigration.
Throughout the paper, for convenience, we always assume
\begin{equation}\label{pa}
  f(0)=g(0)=0,
\end{equation}
so that no extinction happens and for $n\ge 1,$ there are at least $n$ individuals in generation $n.$

 Let $\mathbb T$ be a typical tree generated by the above branching process with immigration.
 For $2\le i\le n,$  pick $i$ individuals from the $n$-th generation at random without replacement and trace their lines of descent back in time till the generation,  marked as $X_{i,1}^n,$ at which those $i$ individuals coalesce into 1 individual for the first time. We call the generation number $X_{i,1}^n$ the coalescence time. We emphasize that due to the existence of immigration, those $i$ individuals may have no common ancestor. In this case, we set $X_{i,1}^n=\infty.$
 Our concern is to study the exact distribution of the coalescence time $X_{i,1}^n$ and find its limit as $n\rto.$

 We remark that without (\ref{pa}), it is enough to consider the distribution of $X_{i,1}^n,$ conditioned on the event $\{N_n\ge i\}$. For $n\ge 0,$ denote by
\begin{align}\label{gf}
\varphi_{\overline n}(z):=E(z^{{N}_n}), z\in[0,1]
\end{align}
  the p.g.f. of ${N}_n$.
   Throughout, for a function $h(z)$  we denote by $h^{(j)}(z)$ its $j$-th derivative and for $n\ge1,$ we set $h_n(z)=h(h_{n-1}(z))$ with $h_0(z)=z.$ In (\ref{gf}),  a subscript $\overline n$ is used to emphasize that  $\varphi_{\overline n}(z)$  is not the $n$-fold composition of a function $\varphi(z).$
   Now we state the main results.
\begin{theorem}\label{bcd} For $n\ge i\ge 2,$ $n>m\ge 0,$ we have
\begin{align}\label{mc}
P(X_{i,1}^n=\infty)=1-\sum_{l=1}^{n}\frac{1}{\Gamma(i)}\int_0^1(1-z)^{i-1}f_{l}^{(i)}(z)g'(f_{l}(z))\prod_{s\neq l,s=1}^{n}g(f_{s}(z))dz,
\end{align}
and
\begin{align}\label{za}
&P(m\le X_{i,1}^n<\infty)=\frac{1}{\Gamma(i)}\int_0^1(1-z)^{i-1}f_{n-m}^{(i)}(z)\varphi_{\overline m}'(f_{n-m}(z))\prod_{l=1}^{n-m}g(f_{l}(z))dz \nonumber\\
&\quad +\sum_{k=1}^{n-m}\frac{1}{\Gamma{(i)}}\int_0^1(1-z)^{i-1}f_{k}^{(i)}(z)g'(f_{k}(z))\prod_{l\neq k,l=1}^{n}g(f_{l}(z))dz.
\end{align}
\end{theorem}

\begin{remark}

 (i) What we have seen from Theorem \ref{bcd}  is the distribution of the coalescence time $X_{i,1}^n$ of $i$ individuals into $1$ individual of the Galton-Watson process with immigration. Basically speaking, for the proof of Theorem \ref{bcd}, we adopt a method similar to \cite{gh}. Of course, because the existence of the immigration, things become more complicated and some difference arises. For example, if those $i$ individuals are chosen from the descendants of  different immigrants, then they will never coalesce into 1 individual. So special light should be shed on the case $X_{i,1}^n=\infty.$ (ii)  We have  results for the coalescence time $X_{i,j}^n,$ that is, the coalescence of $i$ individuals into $j(<i)$ individuals of the Galton-Watson process with immigration. The proofs are more or less similar to those of $X_{i,1}^n$ but the notations are very heavy, so we omit that part in this paper.
\end{remark}

With the distribution of $X_{i,1}^n$ in (\ref{mc}) and (\ref{za}) in hands, one may ask whether $X_{i,1}^n$ converges to a certain distribution as $n\rto.$  This is our next concern. We concentrate only on the limit distribution of  $X_{2,1}^n$ since that of $X_{i,1}^n$ can be treated  similarly. To begin with, we introduce some useful notations. For $m<n,$  let
\begin{align*}
  N_{m,n}^{(l)}:=&\text{the number of progenies in the }n\text{-th generation }\text{born to the }l\text{-th}\nonumber\\
  & \text{individual in the }m\text{-th generation;}\\
  \overline N_{m,n}^{(l)}:=&\text{the number of progenies in the }n\text{-th generation}\text{ born to the }l\text{-th}\nonumber\\
  & \text{immigrant  at generation  }m.
\end{align*}

\begin{lemma}\label{zsam} Assume $f(0)=g(0)=0,$ $g'(1)<\infty,$ $f''(1)<\infty$ and write $\mu=f'(1).$ Then there exist random variables $W,$ $X$ and $V,$ which satisfy $P$-a.s., $0<W,X,V<\infty,$  such that $P$-a.s.,
\begin{align}\label{tc}
 \frac{1}{\mu^n}\overline N_{0,n}^{(1)}\rightarrow W,\  \frac{1}{\mu^{n}}\sum_{k=0}^{n-1}\sum_{l=1}^{I_{k}}\overline N_{k,n}^{(l)}\rightarrow X, \ \frac{1}{\mu^{2n}}\sum_{k=0}^{n-1}\sum_{l=1}^{I_{k}}\z(\overline N_{k,n}^{(l)}\y)^{2}\rightarrow V,
\end{align}
as $n\rto.$ Furthermore, we have $P(V<X^2)=1.$
\end{lemma}
In the next theorem, we study both the quenched and annealed limit laws of the coalescence time $X_{2,1}^n.$

\begin{theorem}\label{sdf} Assume $f(0)=g(0)=0,$ $g'(1)<\infty,$ $f''(1)<\infty.$
Then

{\rm(i)} for almost all family tree $\mathbb T,$ we have
\begin{align}
\lim_{n\to\infty }P(m\le X_{2,1}^n<\infty|\mathbb T)=\frac{\sum_{l=1}^{N_{m}}W_{m,l}^2+V_{m}}{\z(\sum_{l=1}^{N_{m}}W_{m,l}+X_{m}\y)^2},\label{kl}
\end{align}
for $0\leq m<\infty,$  where for $m,l\ge0,$ $W_{m,l},$ $V_m,$ and $X_m$ are independent copies of $W,$ $V$ and $X,$ which appear in Lemma \ref{zsam}, respectively;

{\rm(ii)} for $0\le m<\infty,$
\begin{align}
\lim_{n\to\infty }P(m\le X_{2,1}^n<\infty)=E\z(\frac{\sum_{l=1}^{N_{m}}W_{m,l}^2+V_{m}}{\z(\sum_{l=1}^{N_{m}}W_{m,l}+X_{m}\y)^2}\y),\label{akl}
\end{align}
and especially, \begin{align}
\lim_{n\to\infty }P(X_{2,1}^n<\infty)=E\z(\frac{V}{X^2}\y)<1.\label{mkl}
\end{align}
\end{theorem}
\begin{remark}
  (i) We use a method similar to \cite{a12b} to prove Theorem \ref{sdf}. But, we require $f''(1)<\infty$ or equivalently, $\sigma^2:=\rm{var}(M)<\infty,$ which is not required in \cite{a12b} when studying the coalescence time of Galton-Watson processes without immigration. In our case, by Lemma \ref{wlct} below, $X_{2,1}^n\ge m$ if and only if those two individuals are chosen from progenies born to one of the $N_m$ individuals in the $m$-th generation or from progenies born to the immigrants which immigrate into the system between the $m$-th and $(n-1)$-th generations. Therefore, we need to show the convergence of $ \frac{1}{\mu^{2(n-m)}}\sum_{k=m}^{n-1}\sum_{l=1}^{I_{k}}\z(\overline N_{k,n}^{(l)}\y)^{2},$ which is a submartingale. To apply the convergence theorem of submartingale, we should require the finiteness of $\sigma^2.$ (ii) We see from (\ref{mkl}) that $\lim_{n\to\infty }P(X_{2,1}^n<\infty)<1.$ This happens naturally since for any $n,$ with positive probability, $X_{2,1}^n=\infty,$ that is, those two individuals chosen from progenies in generation $n$ may never coalesce into $1$ individual.
\end{remark}

The remainder of the paper is organized as follows. In Section \ref{eg}, we give two computable examples. Then in Section \ref{pr}, we  provide some lemmas which are useful to prove the main result. Also Lemma \ref{zsam} is proved in this section. Finally,  we devote Section \ref{sep} to finishing the proofs of Theorem \ref{bcd} and Theorem \ref{sdf}.
\section{Examples}\label{eg}
In this section, we give two computable examples.  For two numbers $n\ge k\ge 1,$ we denote by
$$(n)_k=n(n - 1)\cdots(n-k+1),$$ the falling factorial of $n.$
\begin{example}[$l$-nary tree with $k$ immigrants in each generation] \label{ex1}Let $g(z)=z^{k},k\ge 1$ and $f(z)=z^{l},l\ge 2.$ Then by some easy computation, it follows from (\ref{mc}) and (\ref{za}) that for $n>0,$
\begin{align*}
&P(X_{i,1}^n=\infty)=1-\sum_{t=1}^{n}\sum_{s=0}^{i-1}\frac{(-1)^{i-s-1}\tbinom{i-1}{s}(l^t)_ik(1-l)}{(i-1)!(kl(1-l^{n})-s(1-l))},\\
&P(m\le X_{i,1}^n<\infty)=\sum_{s=0}^{i-1}\frac{(-1)^{i-s-1}\binom{i-1}{s}(l^{n-m})_ikl(1-l^m)}{(i-1)!(kl(1-l^{n})-s(1-l))}\\
&\quad\quad\quad\quad\quad\quad\quad\quad\quad\quad+\sum_{t=1}^{n-m}\sum_{s=0}^{i-1}\frac{(-1)^{i-s-1}\tbinom{i-1}{s}(l^t)_ik(1-l)}{(i-1)!(kl(1-l^{n})-s(1-l))}.
\end{align*}
Especially, setting  $i=2,$ then
\begin{align*}
&P(X_{2,1}^n=\infty)=1-\frac{k(1-l)\left(l^2(1-l^{2n})-l(l+1)(1-l^{n})\right)}{\left(1+l\right)
\left(kl(1-l^{n})-(1-l)\right)\left(kl(1-l^{n})\right)},\\
&P(m\le X_{2,1}^n<\infty)=\frac{l^{n-m}(l^{n-m}-1)(1-l^{m})(1-l)}{\left(kl-kl^{n+1}-(1-l)\right)\left(1-l^{n}\right)}\\
&\quad\quad\quad\quad\quad\quad\quad\quad+\frac{k(1-l)\left(l^2(1-l^{2(n-m)})-l(l+1)(1-l^{n-m})\right)}{\left(1+l\right)\left(kl(1-l^{n})-(1-l)\right)\left(kl(1-l^{n})\right)},
\end{align*}
and consequently,
\begin{align*}
&\lim_{n\to\infty }P(X_{2,1}^n=\infty)=\frac{(1+l)k-(l-1)}{(1+l)k},\\
&\lim_{n\rto}P(m\le X_{2,1}^n<\infty)=l^{-2m}(l-1)\frac{1}{k}\z(\frac{l^m-1}{l}+\frac{1}{1+l}\y).
\end{align*}\qed
\end{example}

We see from Example \ref{ex1} that due to the existence of immigration, the formulae for an $l$-nary deterministic tree with $k$ immigrants in each generation are already very complicated. In general, it is difficult to give explicit formulae  if the involved offsprings and immigrations are random. In the next example, we consider  a binary tree with certain random immigration.

\begin{example}
Let $g(z)=\frac{1}{2}(z^2+z)$ and $f(z)=z^2.$ Then some tedious computations  from (\ref{mc}) and (\ref{za}) yield that
\begin{align}\label{ug}
P(m &\le X_{i,1}^{n}<\infty)\no\\
&=\sum_{s=0}^{i-1}\sum_{l=1}^m\sum_{j=0}^{2^{n-m}-1}\sum_{h=0}^{2^{l-1}-1}\sum_{k=0}^{2^{n-l}-1}
\frac{\binom{i-1}{s}(-1)^{i-1-s}(2^{n-m})_i2^{-n+l+1}}{(i-1)!(A(n,m,s,l,j,h,k)+2^{n-m+l})}\no\\
&\quad+\sum_{s=0}^{i-1}\sum_{l=1}^m\sum_{j=0}^{2^{n-m}-1}\sum_{h=0}^{2^{l-1}-1}\sum_{k=0}^{2^{n-l}-1}
\frac{\binom{i-1}{s}(-1)^{i-1-s}(2^{n-m})_i2^{-n+l}}{(i-1)!A(n,m,s,l,j,h,k)}\no\\
&\quad+\sum_{k=1}^{n-m}\sum_{s=0}^{i-1}\sum_{j=0}^{2^{k-1}-1}\sum_{l=0}^{2^{n-k}-1}
\frac{\binom{i-1}{s}(-1)^{i-1-s}(2^{k})_i2^{-n+1}}{(i-1)!\big(B(n,s,l,j,k)+2^k\big)}\no\\
&\quad+\sum_{k=1}^{n-m}\sum_{s=0}^{i-1}\sum_{j=0}^{2^{k-1}-1}\sum_{l=0}^{2^{n-k}-1}
\frac{\binom{i-1}{s}(-1)^{i-1-s}(2^{k})_i2^{-n}}{(i-1)!B(n,s,l,j,k)}
\end{align}
where we write  $A(n,m,s,l,j,h,k):=2^{n+1}-2-s+2j+h2^{n-m+1}+k2^{n-m+l+1}$ and  $B(n,s,l,j,k):=2^{n+1}+l2^{k+1}+2j-2-s$ for simplicity. \qed

We remark that though the formula in (\ref{ug}) looks very ugly, given $n,m$ and $i,$ it is indeed computable.
\end{example}

\section{Preliminary results}\label{pr}

In this section, we present several lemmas which are useful to prove  Theorem \ref{bcd} and \ref{sdf}. To begin with, we study the quenched distribution of $X_{i,1}^n.$
\begin{lemma}\label{wlct}
We have
\begin{align}\label{qpm}
P\z(m\le X_{i,1}^n<\infty|\mathbb T\y)=\frac{\sum_{l=1}^{N_{m}}\z(N_{m,n}^{(l)}\y)_{i}+\sum_{k=m}^{n-1}\sum_{l=1}^{I_k}\z(\overline N_{k,n}^{(l)}\y)_i}{(N_n)_i},
\end{align}
and especially,
\begin{align}\label{mpi}
P\z(X_{i,1}^n<\infty|\mathbb T\y)=\frac{\sum_{k=0}^{n-1}\sum_{l=1}^{I_k}\z(\overline N_{k,n}^{(l)}\y)_i}{(N_n)_i}.
\end{align}
\end{lemma}
\proof Note that the event $\{m\le X_{i,1}^n<\infty\}$ occurs if and only if those $i$ individuals are randomly chosen from progenies in the $n$-th generation born to one individual in the $m$-th generation or born to one of those immigrants which immigrate into the system between the $m$-th and $(n-1)$-th generation. Thus (\ref{qpm}) follows.  Since $N_0=0,$ letting $m=0$ in (\ref{qpm}), we get (\ref{mpi}).\qed

To obtain the annealed distribution of $X_{i,1}^n,$ it is necessary to take expectation in both sides of (\ref{qpm}) and (\ref{mpi}). The following lemma supplies a tool to compute those expectations.
\begin{lemma}\label{zaa}
Suppose that $Y_1,...,Y_n, Z$ are mutually independent $\mathbb Z^+$-valued random variables, $Y_1,...,Y_n$ share the same p.g.f. $\phi(z),$ and $Z$ has p.g.f. $h(z).$
Then, for $1 \le j< i,$ $i_1+i_2+...+i_j=i$ and $i_1,\;i_2,...\;i_j\geq 1,$
\begin{align*}
E\left(\frac{(Y_1)_{i_1}(Y_2)_{i_2}\cdots(Y_j)_{i_j}}{(\sum_{l=1}^{n}Y_{l}+Z)_i}\right)
=\frac{1}{\Gamma{(i)}}\int_0^1(1-z)^{i-1}\phi(z)^{n-j}\prod_{l=1}^{j}\phi^{(i_l)}(z)h(z)dz.
\end{align*}

\end{lemma}
\proof Setting  $\sum_{l=j+1}^{n}Y_{l}=\mathcal{Y},$ we have
\begin{align}
&E\left(\frac{(Y_1)_{i_1}(Y_2)_{i_2}\cdots(Y_j)_{i_j}}{(\sum_{l=1}^{n}Y_{l}+Z)_i}\right)\no\\
&=\sum_{m,m_1,...,m_j,x=0}^{\infty}P(\mathcal{Y}=m)P(Y_1=m_1)\cdots P({Y}_j=m_j)P(Z=x)\no\\
&\quad \times\frac{(m_1)_{i_1}(m_2)_{i_2}\cdots(m_j)_{i_j}}{(m+m_1+...+m_j+x)_i}.\label{qa}
\end{align}
On the other hand, owing to the independence of $Y_1,...,Y_n, Z,$ some careful computation yields that
\begin{align}
&\frac{1}{\Gamma{(i)}}\int_0^1(1-z)^{i-1}\phi(z)^{n-j}\prod_{l=1}^{j}\phi^{(i_l)}(z)h(z)dz\no\\
&=\sum_{m,m_1,...,m_j,x=0}^{\infty}P(\mathcal{Y}=m)P({Y}_1=m_1)\cdots P({Y}_j=m_j)P(Z=x)\no\\
&\quad \quad \quad \quad \quad \times(m_1)_{i_1}(m_2)_{i_2}...(m_j)_{i_j}\frac{1}{\Gamma{(i)}}\int_0^1(1-z)^{i-1}z^{m+m_1+...+m_j+x-i}dz\no\\
&=\sum_{m,m_1,...,m_j,x=0}^{\infty}P(\mathcal{Y}=m)P({Y}_1=m_1)\cdots P({Y}_j=m_j)P(Z=x)\no\\
&\quad \times\frac{(m_1)_{i_1}(m_2)_{i_2}\cdots (m_j)_{i_j}}{(m+m_1+...+m_j+x)_i}.\label{qb}
\end{align}
In view of (\ref{qa}) and (\ref{qb}), Lemma \ref{zaa} follows.\qed

At last, we give the proof of Lemma \ref{zsam} to conclude this section.

\vspace{0.5cm}

\noindent {\it Proof of Lemma \ref{zsam}.}

Write $\lambda:=g'(1),$ $\sigma^2:=f''(1)+f'(1)-f'(1)^2\equiv\mathrm{Var}(M)$ and for $n\ge1$ set   $$W_n:=\frac{\overline N_{0,n}^{(1)}}{\mu^n}, \ X_n:=\frac{1}{\mu^{n}}\sum_{k=0}^{n-1}\sum_{l=1}^{I_{k}}\overline N_{k,n}^{(l)},\ V_n:= \frac{1}{\mu^{2n}}\sum_{k=0}^{n-1}\sum_{l=1}^{I_{k}}\z(\overline N_{k,n}^{(l)}\y)^{2}.$$
By \cite{an}, we have $P$-a.s., $W_n\rightarrow W$ as $n\rto$ for some  random variable $W.$   Since we assume $f(0)=0,$ no extinction happens to the branching process $\overline N_{0,n}^{(1)},n\ge 1$ so that $P(0<W<\infty)=1.$

 The convergence of $X_n$ and $V_n$ can be shown by the convergence theorem of submartingale. We remark that the convergence of $X_n$ is a direct consequence of \cite{se}.
We give its proof here for convenience. Some direct computation yields that
\begin{align*}
  &E(X_n)=\frac{\mu(\mu^n-1)}{\mu^n(\mu-1)},\\
  &E(V_{n})=\frac{\lambda}{\mu^{2n}}\Big(\frac{\sigma^2}{\mu-1}\Big(\frac{\mu(\mu^{2n}-1)}{\mu^2-1}-\frac{\mu^{n}-1}{\mu-1}\Big)+\frac{\mu^2(\mu^{2n}-1)}{\mu^2-1}\Big).
\end{align*}
Thus there exists some constant $C>0,$ such that for all $n\ge1,$
\begin{align}\label{ecb}E(X_n)<C,\ E(V_n)<C.\end{align}
Notice that by Markov property, \begin{align}\label{np}
&E(V_{n+1}|V_{n},V_{n-1},...,V_{1})=E(V_{n+1}|V_{n})\no\\
&=\frac{1}{\mu^{2n+2}}E\Big(\sum_{k=0}^{n-1}\sum_{l=1}^{I_{k}}\z(\overline N_{k,n+1}^{(l)}\y)^{2}+\sum_{l=1}^{I_{n}}\z(\overline N_{n,n+1}^{(l)}\y)^{2}\Big|V_{n}\Big)\no\\
&=\frac{1}{\mu^{2n+2}}\Big[E\Big(\sum_{l=1}^{I_{n}}\z(\overline N_{n,n+1}^{(l)}\y)^{2}\Big)+E\Big(\sum_{k=0}^{n-1}\sum_{l=1}^{I_{k}}\Big(\sum_{j=1}^{\overline N_{k,n}^{(l)}}\xi_{n,k,l,j}\Big)^2\Big|V_{n}\Big)\Big],
\end{align}
where $\xi_{n,k,l,j}, n,k,l,j\ge1$ are mutually independent random variables whose p.g.f. are all $f(z).$ Since
$E\z(\sum_{l=1}^{I_{n}}\z(\overline N_{n,n+1}^{(l)}\y)^{2}\y)=\lambda(\sigma^2+\mu^2)$
and \begin{align*}
E\Big(\sum_{k=0}^{n-1}&\sum_{l=1}^{I_{k}}\Big(\sum_{j=1}^{\overline N_{k,n}^{(l)}}\xi_{n,k,l,j}\Big)^2\Big|V_{n}\Big)\no\\
&=\sum_{k=0}^{n-1}\sum_{l=1}^{I_{k}}\z(\sum_{j=1}^{\overline N_{k,n}^{(l)}}E(\xi_{n,k,l,j}^2)+\sum_{i,j=1,i\neq j}^{\overline N_{k,n}^{(l)}}E(\xi_{n,k,l,i})E(\xi_{n,k,l,j})\y)\no\\
&=\sum_{k=0}^{n-1}\sum_{l=1}^{I_{k}}\z(\z(\overline N_{k,n}^{(l)}\y)^2\mu^2+\overline N_{k,n}^{(l)}\sigma^2\y),
\end{align*}
from (\ref{np}), we get
\begin{align*}
&E(V_{n+1}|V_{n},V_{n-1},...,V_{1})\no\\
&=\frac{1}{\mu^{2n+2}}\z(\lambda(\sigma^2+\mu^2)+\sum_{k=0}^{n-1}\sum_{l=1}^{I_{k}}\z(\z(\overline N_{k,n}^{(l)}\y)^2\mu^2+\overline N_{k,n}^{(l)}\sigma^2\y)\y)\no\\
&>\frac{1}{\mu^{2n}}\sum_{k=0}^{n-1}\sum_{l=1}^{I_{k}}\z(\overline N_{k,n}^{(l)}\y)^2=V_{n}.
\end{align*}
A similar argument also yields that $E(X_{n+1}|X_n,...,X_1)>X_n.$
Consequently, taking (\ref{ecb}) into consideration, we conclude that both $V_{n},$ $n\geq1$ and $X_n,n\ge1$ are submartingales with bounded means. Thus,
applying the convergence theorem of submartingale, we have $P$-a.s., $$\lim_{n\to\infty }V_{n}=V, \ \lim_{n\rto}X_n=X $$ for some  random variables $V$ and $X$ with $E(V)<C$ and $ E(X)<C.$
Clearly, we have for all $n\ge1$ $X_n>\frac{\overline N_{0,n}^{(1)}}{\mu^n}$ and $V_n>\frac{\overline N_{0,n}^{(1)}}{\mu^n}.$ Therefore, we have $P$-a.s., $0<V,X<\infty.$
We thus finish the proof of (\ref{tc}).
Finally,  by (\ref{tc}), we have $P$-a.s.,
\begin{align*}
  X^2-V=\lim_{n\rto}(X_n^2-V_n)\ge\lim_{n\rto}\frac{ N_{0,n}^{(1)}}{\mu^n}\frac{ N_{0,n}^{(2)}}{\mu^n}=W^2>0.
\end{align*}
 Lemma \ref{zsam} is proved. \qed

\section{Proofs of main results}\label{sep}
\subsection{Proof of Theorem \ref{bcd}}
To begin with, we show (\ref{za}). By Lemma \ref{wlct}, we have
\begin{align}\label{mif}
P\z(m\le X_{i,1}^n<\infty\y)&=E\z(\frac{\sum_{l=1}^{N_{m}}\z(N_{m,n}^{(l)}\y)_{i}+\sum_{t=m}^{n-1}\sum_{l=1}^{I_t}\z(\overline N_{t,n}^{(l)}\y)_i}{(N_n)_i}\y)\no\\
&=E\z(\frac{\sum_{l=1}^{N_{m}}\z(N_{m,n}^{(l)}\y)_{i}}{(N_n)_i}\y)+\sum_{t=m}^{n-1}E\z(\frac{\sum_{l=1}^{I_t}\z(\overline N_{t,n}^{(l)}\y)_i}{(N_n)_i}\y)\no\\
&=:\mathcal{L}+\sum_{t=m}^{n-1}\mathcal{L}_t.
\end{align}

To compute $\mathcal L,$ noticing that ${N_n}=\sum_{l=1}^{N_m}N_{m,n}^{(l)}+\sum_{k=m}^{n-1}\sum_{l=1}^{I_{k}}{\overline N_{{k},n}^{(l)}},$ thus,
\begin{align}\label{ca}
  \mathcal{L}&=E\z(\frac{\sum_{l=1}^{N_{m}}\z(N_{m,n}^{(l)}\y)_{i}}{(N_n)_i}\y)\no\\
  &=\sum_{n_m,i_m,...,i_{n-1}=0}^{\infty}P(N_m=n_m)\prod_{k=m}^{n-1}P(I_k=i_k)\nonumber\\
&\quad\quad\quad\quad \times \sum_{l=1}^{n_m}E\left(\frac{(N_{m,n}^{(l)})_i}{\left(\sum_{l=1}^{n_m}N_{m,n}^{(l)}+\sum_{k=m}^{n-1}\sum_{l=1}^{i_{k}}{\overline N_{{k},n}^{(l)}}\right)_{i}}\right).
\end{align}
In order to apply Lemma \ref{zaa}, we set $Z=\sum_{k=m}^{n-1}\sum_{l=1}^{i_{k}}{\overline N_{{k},n}^{(l)}}$ and $Y_l=N_{m,n}^{(l)}, 1\le l\le n_m.$ Then
$Z,Y_1,...,Y_{n_m}$ are mutually independent. By some classical argument of Galton-Watson processes, $Y_1,...,Y_{n_m}$ share the common p.g.f. $f_{n-m}(z), z\in[0,1].$  From the assumption, it is clear that  $N_{{k},n}^{(l)}, m\le k\le n-1,1\le l\le i_k$ are mutually independent and for $m\le k\le n-1,1\le l\le i_k,$ the p.g.f. of $N_{{k},n}^{(l)}$  equals $f_{n-k}(z),z\in[0,1].$ Thus, the p.g.f. of $Z$ equals $\prod_{k=m}^{n-1} f_{n-k}^{i_k}(z),z\in[0,1].$ Therefore, an application of
Lemma \ref{zaa} yields that
\begin{align}\label{cc}
&E\left(\frac{(N_{m,n}^{(l)})_i}{\left(\sum_{l=1}^{n_m}N_{m,n}^{(l)}+\sum_{k=m}^{n-1}\sum_{l=1}^{i_{k}}{\overline N_{{k},n}^{(l)}}\right)_{i}}\right)\nonumber\\
&=\frac{1}{\Gamma{(i)}}\int_0^1(1-z)^{i-1}f_{n-m}^{(i)}(z)f_{n-m}^{n_m-1}(z)\prod_{k=m}^{n-1}f_{n-k}^{i_{k}}(z)dz.
\end{align}
Substituting (\ref{cc}) into (\ref{ca}), we get that
\begin{align}\label{cd}
\mathcal{L}&=\sum_{n_m,i_m,...,i_{n-1}=0}^{\infty}P(N_m=n_m)\prod_{k=m}^{n-1}P(I_k=i_k)\nonumber\\
&\quad\quad\times \frac{n_m}{\Gamma{(i)}}\int_0^1(1-z)^{i-1}f_{n-m}^{(i)}(z)f_{n-m}^{n_m-1}(z)\prod_{k=m}^{n-1}f_{n-k}^{i_{k}}(z)dz \nonumber\\
&=\frac{1}{\Gamma{(i)}}\int_0^1(1-z)^{i-1}f_{n-m}^{(i)}(z)\varphi_{\overline m}'(f_{n-m}(z))
\prod_{k=m}^{n-1}g(f_{n-k}(z))dz\no\\
&=\frac{1}{\Gamma{(i)}}\int_0^1(1-z)^{i-1}f_{n-m}^{(i)}(z)\varphi_{\overline m}'(f_{n-m}(z))
\prod_{k=1}^{n-m}g(f_{k}(z))dz.
\end{align}
$\mathcal L_t$ can be computed similarly. Indeed,
\begin{align}\label{ze}
\mathcal L_t&=E\z(\frac{\sum_{l=1}^{I_t}\z(\overline N_{t,n}^{(l)}\y)_i}{(N_n)_i}\y)=E\left(\frac{\sum_{l=1}^{I_t}\z(\overline N_{t,n}^{(l)}\y)_i}{\z(\sum_{l=1}^{N_t}N_{t,n}^{(l)}+\sum_{k=t}^{n-1}\sum_{l=1}^{I_{k}}{\overline N_{{k},n}^{(l)}}\y)_i}\right) \nonumber\\
&=\sum_{n_t,i_t,...,i_{n-1}=0}^{\infty}P(N_t=n_t)\prod_{k=t}^{n-1}P(I_k=i_k)\no\\
&\quad \quad \quad \quad \times  \sum_{l=1}^{i_t}E\left(\frac{\z(\overline N_{t,n}^{(l)}\y)_i}{\z(\sum_{l=1}^{n_t}N_{t,n}^{(l)}+\sum_{k=t}^{n-1}\sum_{l=1}^{i_{k}}{\overline N_{{k},n}^{(l)}}\y)_i}\right).
\end{align}
Now setting $Y_l=\overline N_{t,n}^{(l)}, 1\le l\le i_t$ and $Z=\sum_{l=1}^{n_t}N_{t,n}^{(l)}+\sum_{k=t+1}^{n-1}\sum_{l=1}^{i_{k}}{\overline N_{{k},n}^{(l)}},$ then some similar arguments as above yield that $Z,Y_1,..,Y_{i_t}$ are mutually independent, $Y_1,..,Y_{i_t}$ share the common p.g.f. $f_{n-t}(z), z\in[0,1],$ and the p.g.f. of $Z$ equals $f_{n-t}^{n_t}(z)\prod_{k=t+1}^{n-1}f_{n-k}^{i_{k}}(z), z\in[0,1].$
Applying Lemma \ref{zaa} again, we have
\begin{align}\label{abh}
&E\z(\frac{\z(\overline N_{t,n}^{(l)}\y)_i}{\z(\sum_{l=1}^{n_t}N_{t,n}^{(l)}+\sum_{k=t}^{n-1}\sum_{l=1}^{i_{k}}{\overline N_{{k},n}^{(l)}}\y)_i}\y)\no\\
&\quad\quad=\frac{1}{\Gamma{(i)}}\int_0^1(1-z)^{i-1}f_{n-t}^{(i)}(z)f_{n-t}^{i_t+n_t-1}(z)\prod_{k=t+1}^{n-1}f_{n-k}^{i_{k}}(z)dz.
\end{align}
Substituting (\ref{abh}) into (\ref{ze}), we get
\begin{align}\label{eb}
\mathcal{L}_t&=\sum_{n_t,i_t,...,i_{n-1}=0}^{\infty}P(N_t=n_t)\prod_{k=t}^{n-1}P(I_k=i_k)\no\\
&\quad\quad\quad\quad\quad\quad \times\frac{ i_t}{\Gamma{(i)}}\int_0^1(1-z)^{i-1}f_{n-t}^{(i)}(z)f_{n-t}^{i_t+n_t-1}(z)\prod_{k=t+1}^{n-1}f_{n-k}^{i_{k}}(z)dz\no\\
&=\frac{1}{\Gamma{(i)}}\int_0^1(1-z)^{i-1}f_{n-t}^{(i)}(z)g'(f_{n-t}(z))\varphi_{\overline t}(f_{n-t}(z))\prod_{k=t+1}^{n-1}g(f_{n-k}(z))dz.
\end{align}
But it follows from (\ref{nr}) and (\ref{gf}) that for $n\ge1,$
\begin{align}\label{nn}
\varphi_{\overline n}(z)&=E(z^{N_n})=\varphi_{\overline{n-1}}(f(z))g(f(z))
\end{align}
with $\varphi_{\overline 0}(z)=1.$ Iterating (\ref{nn}), we get
\begin{align*}
\varphi_{\overline n}(z)=\prod_{l=1}^{n}g(f_l(z)),\ n\ge1.
\end{align*}
As a consequence, we obtain
\begin{align}\label{abf}
\varphi_{\overline t}(f_{n-t}(z))=\prod_{l=n-t+1}^{n}g(f_l(z)).
\end{align}
Thus, it follows from (\ref{eb}) and (\ref{abf}) that
\begin{align}\label{zg}
\mathcal L_t&=\frac{1}{\Gamma{(i)}}\int_0^1(1-z)^{i-1}f_{n-t}^{(i)}(z)g'(f_{n-t}(z))\prod_{l\neq{n-t},l=1}^{n}g(f_{l}(z))dz.
\end{align}
Consequently, taking (\ref{mif}), (\ref{cd}) and (\ref{zg}) together, we finish the proof of (\ref{za}).

Next we proceed to prove (\ref{mc}). Letting $m=0$ in (\ref{za}), we get
\begin{align*}
P(X_{i,1}^n<\infty)=\sum_{l=1}^{n}\frac{1}{\Gamma(i)}\int_0^1(1-z)^{i-1}f_{l}^{(i)}(z)g'(f_{l}(z))\prod_{s\neq l,s=1}^{n}g(f_{s}(z))dz.
\end{align*}
Then (\ref{mc}) is proved and so is  Theorem \ref{bcd}.
\qed

\subsection{Proof of Theorem \ref{sdf}}
Fix $0\leq m<\infty.$ Then it follows from Lemma \ref{wlct} that
\begin{align}\label{qpdd}
P(&m\le X_{2,1}^n<\infty|\mathbb T)=\frac{\sum_{l=1}^{N_{m}}N_{m,n}^{(l)}\big(N_{m,n}^{(l)}-1\big)+\sum_{k=m}^{n-1}\sum_{l=1}^{I_{k}}\overline N_{k,n}^{(l)}\big(\overline N_{k,n}^{(l)}-1\big)}{N_n(N_n-1)}\no\\
&=\frac{\sum_{l=1}^{N_{m}}N_{m,n}^{(l)}\big(N_{m,n}^{(l)}-1\big)+\sum_{k=m}^{n-1}\sum_{l=1}^{I_{k}}\overline N_{k,n}^{(l)}\big(\overline N_{k,n}^{(l)}-1\big)}{\big(\sum_{l=1}^{N_{m}}N_{m,n}^{(l)}+\sum_{k=m}^{n-1}\sum_{l=1}^{I_{k}}\overline N_{k,n}^{(l)}\big)\big(\sum_{l=1}^{N_{m}}N_{m,n}^{(l)}+\sum_{k=m}^{n-1}\sum_{l=1}^{I_{k}}\overline N_{k,n}^{(l)}-1\big)}\no\\
&=\Delta_{1}+\Delta_{2},
\end{align}
where \begin{align*}
  &\Delta_{1}:=\frac{\sum_{l=1}^{N_{m}}N_{m,n}^{(l)}(N_{m,n}^{(l)}-1)}{(\sum_{l=1}^{N_{m}}N_{m,n}^{(l)}+\sum_{k=m}^{n-1}\sum_{l=1}^{I_{k}}\overline N_{k,n}^{(l)})(\sum_{l=1}^{N_{m}}N_{m,n}^{(l)}+\sum_{k=m}^{n-1}\sum_{l=1}^{I_{k}}\overline N_{k,n}^{(l)}-1)},\\
  &\Delta_2:=\frac{\sum_{k=m}^{n-1}\sum_{l=1}^{I_{k}}\overline N_{k,n}^{(l)}(\overline N_{k,n}^{(l)}-1)}{(\sum_{l=1}^{N_{m}}N_{m,n}^{(l)}+\sum_{k=m}^{n-1}\sum_{l=1}^{I_{k}}\overline N_{k,n}^{(l)})(\sum_{l=1}^{N_{m}}N_{m,n}^{(l)}+\sum_{k=m}^{n-1}\sum_{l=1}^{I_{k}}\overline N_{k,n}^{(l)}-1)}.
\end{align*}Applying Lemma \ref{zsam}, we get that $P$-a.s.,
\begin{align}
\Delta_{1}&=\frac{\sum_{l=1}^{N_{m}}\frac{N_{m,n}^{(l)}}{\mu^{n-m}}\frac{N_{m,n}^{(l)}-1}{\mu^{n-m}}}{\frac{(\sum_{l=1}^{N_{m}}N_{m,n}^{(l)}+\sum_{k=m}^{n-1}\sum_{l=1}^{I_{k}}\overline N_{k,n}^{(l)})}{\mu^{n-m}}\times\frac{(\sum_{l=1}^{N_{m}}N_{m,n}^{(l)}+\sum_{k=m}^{n-1}\sum_{l=1}^{I_{k}}\overline N_{k,n}^{(l)}-1)}{\mu^{n-m}}}\no\\
&\rightarrow\frac{\sum_{l=1}^{N_{m}}W_{m,l}^2}{(\sum_{l=1}^{N_{m}}W_{m,l}+X_{m})^2} \quad \text{ as } n\rightarrow\infty,\label{ou}
\end{align}
and
\begin{align}
\Delta_{2}&=\frac{\frac{1}{\mu^{2(n-m)}}\sum_{k=m}^{n-1}\sum_{l=1}^{I_{k}}\overline N_{k,n}^{(l)}\big(\overline N_{k,n}^{(l)}-1\big)}{\frac{\sum_{l=1}^{N_{m}}N_{m,n}^{(l)}+\sum_{k=m}^{n-1}\sum_{l=1}^{I_{k}}\overline N_{k,n}^{(l)}}{\mu^{n-m}}\times\frac{\sum_{l=1}^{N_{m}}N_{m,n}^{(l)}+\sum_{k=m}^{n-1}\sum_{l=1}^{I_{k}}\overline N_{k,n}^{(l)}-1}{\mu^{n-m}}}\no\\
&\rightarrow\frac{V_{m}}{(\sum_{l=1}^{N_{m}}W_{m,l}+X_{m})^2}\quad\text{ as } n\rightarrow\infty,\label{op}
\end{align}
where for $m,l\ge0,$ $W_{m,l},$ $V_m$ and $X_m$ are independent copies of $W,$ $V$ and $X,$ which appear in Lemma \ref{zsam}, respectively.

Substituting (\ref{ou}) and (\ref{op}) into (\ref{qpdd}), we obtain (\ref{kl}).
We thus finish the proof of Part (i) of Theorem \ref{sdf}.

Now we proceed to prove the second part. Taking expectation in both sides of (\ref{kl}) and using the bounded convergence theorem, we get (\ref{akl}). Finally, putting especially $m=0$ in (\ref{akl}), we have \begin{align*}
\lim_{n\to\infty }P(0\le X_{2,1}^n<\infty)=E\z(V/X^2\y).
\end{align*}
But by Lemma \ref{zsam}, we have $P(0<V<X^2)=1,$ implying  $E(V/X^2)<1.$  Consequently, Part (ii) of Theorem \ref{sdf} is proved. \qed

\vspace{0.5Cm}

\noindent{\large{\bf \Large Acknowledgements.}} The financial support of the National
Nature Science Foundation of China (Grant No. 11501008) is gratefully acknowledged.


\begin{thebibliography}{99}
\addtolength{\itemsep}{-0.5em}

\bibitem{a12} K. B. Athreya,  Coalescence in critical and subcritical Galton-Watson branching processes, \textit{J. Appl. Probab.} {\bf49} (2012) 627-638.
\bibitem{a12b}K. B. Athreya,  Coalescence in the recent past in rapidly growing populations, \textit{ Stochastic Process. Appl.} \textbf{122} (2012) 3757-3766.
\bibitem{a15} K. B. Athreya,  Coalescence in branching processes, \textit{ Chapter in Branching processes and their applications, Leture Notes in Statistics} \textbf{219} (2015) 3-22.
     \bibitem{ah} K. B. Athreya and J. Hong,  Coalescence on Supercritical Bellman-Harris Branching Processes, {\it Taiwanese J. Math.} \textbf{22} (2018)  245-261.
\bibitem{an}K. B. Athreya and P. E. Ney,  {\it Branching Processes}(New York: Dover 2004).
\bibitem{se} E. Seneta,  A note on the supercritical Galton-Watson process with immigration, {\it Math. Biosci.} \text{6} (1970)  305-311.

\bibitem{gh} N. Grosjean and T. Huillet,  On the genealogy and coalescence times of Bienaym$\acute{e}$-Galton-Watson branching processes, \textit{Stoch. Models}  \textbf{34} (2018) 1-24.

\bibitem{ha} J. Hong, {\it Coalescence in Bellman-Harris and multi-type branching processes} (Iowa State University Capstones, Graduate Theses $\&$ Dissertations 10103, 2011) (https://lib.dr.iastate.edu/etd/10103).


\bibitem{lama} A. Lambert,  Coalescence times for the branching process, \textit{Adv. in Appl. Probab.} \textbf{35} (2003) 1071-1089.
\bibitem{lp} A. Lambert and L. Popovic,  The coalescent point process of branching trees, \textit{Ann. Appl. Probab.,} \textbf{23}  (2013) 99-144.
 \end{thebibliography}
\end{document}